\documentclass[number,citesort,dvips]{arxbj}
\usepackage{cursive}  

\aid{0}
\volume{16}
\issue{4}
\pubyear{2010}
\firstpage{1294}
\lastpage{1311}
\doi{10.3150/10-BEJ261}

\makeatletter
\newtheorem{lemma}{Lemma}
\newtheorem{proposition}[lemma]{Proposition}
\newtheorem{corollary}[lemma]{Corollary}
\newtheorem{theorem}[lemma]{Theorem}

\newremark{remark}[lemma]{Remark}

\newproclaim{definition}[lemma]{Definition}

\def\ve{\varepsilon}
\def\e{\mathbb{E}}
\def\re{\mathbb{R}}
\def\sign{\operatorname{sign}}
\def\Var{\operatorname{Var}}

\makeatother

\begin{document}
\begin{frontmatter}

\title{Local time and Tanaka formula for a Volterra-type
multifractional Gaussian~process}
\runtitle{Local time and Tanaka formula for a multifractional process}
\begin{aug}
\author[a]{\fnms{Brahim} \snm{Boufoussi}\thanksref{a}\ead[label=e1]{boufoussi@ucam.ac.ma}},
\author[b]{\fnms{Marco} \snm{Dozzi}\thanksref{b,e2}\ead[label=e2,mark]{marco.dozzi@iecn.u-nancy.fr}} \and
\author[b]{\fnms{Renaud} \snm{Marty}\thanksref{b,e3}\ead[label=e3,mark]{renaud.marty@iecn.u-nancy.fr}\corref{}}
\runauthor{B. Boufoussi, M. Dozzi and R. Marty}
\address[a]{Facult\'e des Sciences Semlalia, Department of Mathematics,
Cadi Ayyad University, BP 2390 Marrakech, Morocco. \printead{e1}}
\address[b]{Institut Elie Cartan de Nancy, Nancy-Universit\'e, CNRS,
INRIA, B.P. 239, F-54506 Vandoeuvre-l\`es-Nancy Cedex, France. \printead{e2}; \printead*{e3}}
\end{aug}

\received{\smonth{5} \syear{2009}}
\revised{\smonth{1} \syear{2010}}

%
\begin{abstract}
The stochastic calculus for Gaussian processes is applied to obtain a Tanaka
formula for a Volterra-type multifractional Gaussian process. The existence
and regularity properties of the local time of this process are obtained
by means of Berman's Fourier analytic approach.
\end{abstract}

%

\end{frontmatter}
%

\section{Introduction}

Several types of multifractional Gaussian processes have been studied,
including the
processes introduced by L\'{e}vy-V\'{e}hel and Peltier \cite{LVP} and Benassi,
Jaffard and Roux \cite{BJR}, known as the moving average and harmonizable
versions of multifractional Brownian motion, and the processes
introduced by
Benassi, Cohen and Istas \cite{BCI} and by Surgailis \cite{S}. These
processes are usually defined by replacing, in certain representations
of fractional
Brownian motion (fBm), the Hurst parameter $H$ by a Hurst function $h$,
that is, a real
function of the time parameter with values in (0,1).
Generalizations to stable processes have also been studied \cite{L,ST1,ST2}.
The stochastic properties and the regularity of the trajectories of
these processes
can be characterized in the same terms as for fBm, in particular, by
the property of
local asymptotic self-similarity or by the pointwise H\"older exponent,
which is
constant and equal to $H$ for fBm, but which varies in
time following $h$ in the multifractional case.

The aim of this article is to study a Volterra-type multifractional Gaussian
process $B_h$ which fits into the framework of stochastic calculus, and
its local time.
In fact, the stochastic calculus, including the It\^{o} formula and stochastic
differential equations, is now well established for fBm and
stochastic integrals have been defined in the Malliavin sense or by
means of Wick
products. The stochastic calculus for multifractional Gaussian
processes has not yet
been developed explicitly, but important elements, including a
divergence integral
and an It\^{o} formula, have been proven in \cite{AMN} for a class of Gaussian
processes admitting a kernel representation with respect to Brownian
motion. The
process $B_h$ we study in this article is defined by allowing for a
Hurst function
in the kernel representation of fBm and it belongs to the class of Gaussian
processes studied in \cite{AMN}. The It\^{o} formula in \cite{AMN} is
applied to
obtain a representation of the local time for which, as is usually the
case in
Tanaka-like formulas, the occupation measure is not the Lebesgue
measure, but the
quadratic variation process. However, contrary to fBm, the
interpretation in our case is
more delicate because the quadratic variation process is not
necessarily increasing
(but of bounded variation). We compare this local time to the local
time with
respect to the Lebesgue measure, which we obtain, together with the regularity
properties of its trajectories, by the Fourier analytic approach
initiated by Berman \cite{Ber69} based on the notion of local
non-determinism (LND).

It is well known that the cases where the Hurst parameter
or the Hurst function is $<$1$/$2  (resp., $>$1$/$2) have to
be treated separately (see \cite{AMN} and \cite{N}). We are interested
in the case where
the Hurst function $h$ takes values in (1$/$2,1). In fact, the Volterra-type
process we study here is defined only for this case, which is appropriate
for long memory applications \cite{SamoTaqqu}.
The article is organized as follows. In Section 2, stochastic
properties and the
regularity of the trajectories of~$B_h$ are proved. The form of the covariance
function (Proposition 2) shows, in particular, that~$B_h$ differs from the
harmonizable and moving average multifractional Brownian motions. The
continuity of
the trajectories of~$B_h$ is obtained by classical criteria from
estimates for the
second order moments of the increments (Proposition 3). The lass
property is proved
in Proposition 5; it implies that the pointwise H\"{o}lder exponent of~$B_h$ is
equal to $h$ (Proposition 6). In Section 3, it is shown that $B_h$
satisfies the property
labeled (K4) in \cite{AMN} and that the quadratic variation of~$B_h$ is
of bounded
variation. Therefore, the divergence integral and the It\^{o} formula
developed in
\cite{AMN} (recalled briefly in this section) hold for $B_h$. Section 4
is devoted
to the local time of $B_h$. Its existence and square integrability with
respect to
the space variable follow from a classical criterion due to Berman
\cite{Ber69} (Proposition 11). A Tanaka-type formula with the divergence
integral of Section 3 is then given (Theorem 12), where the occupation
measure of the
local time is the quadratic variation of $B_h$. Finally, the local time
with respect
to the Lebesgue measure and its regularity in the space and time
variables are given
in Theorem 16. The proofs are based on the LND property, which is shown
to hold for
$B_h$ in Proposition 14.

\section{A Volterra-type multifractional Gaussian process}

It is well known that the fractional Brownian motion $B_{H}$ with (fixed)
Hurst parameter $H\in({1/2},1)$ can be represented for any $t\geq0$ as
\[
B_{H}(t) =\int_{0}^{t}K_{H}(t,u)W(\mathrm{d}u),
\]
where
%
\begin{equation}
\label{defKH}
K_{H}(t,u)=u^{
{1/2}-H}\int_{u}^{t}(y-u)^{H-{3/2}}y^{H-{1/2}}\,\mathrm{d}y
\end{equation}
and $W(\mathrm{d}y)$ is a Gaussian measure.
Let $a$ and $b$ be two real numbers satisfying
$1/2<a<b<1$. Throughout the paper, we consider a function
$h\dvtx \mathbb{R} \rightarrow[a,b] $. We assume that this function is
$\beta$-H\"older with
$\sup h< \beta$.
Define the centered Gaussian process
$B_{h}=\{B_{h}(t), t\geq0\}$ by
\[
B_{h}(t)=B_{h(t)}(t)=\int_{0}^{t}K_{h(t)}(t,u)W(\mathrm{d}u),
\]
where
%
\begin{equation}
\label{defKh}
K_{h(t)}(t,u)=u^{{1/2}-h(t)}\int_{u}^{t}(y-u)^{h(t)-{3/2}
}y^{h(t)-{1/2}}\,\mathrm{d}y.
\end{equation}
Before establishing properties of $B_h$, we state a lemma regarding
an estimate on $K_H$ that we use throughout the paper.
\begin{lemma}
\label{estimateKH}
For every $T>0$,
there exists a function $\Phi_T\in L^2((0,T], \re_+)$ such that for
every $s\in(0,T]$,
\[
\sup_{\lambda\in[a,b], t\in(0,T]}\biggl| \frac{
\partial}{\partial\lambda}K_{\lambda}(t,s)\biggr| \leq\Phi_T(s).
\]
\end{lemma}
\begin{pf}
We have
\begin{eqnarray*}
\frac{\partial}{\partial\lambda}K_{\lambda}(t,s) & = & (-\!\log
s)s^{1/2-\lambda
}\int_{s}^{t}(y-s)^{\lambda-3/2}y^{\lambda-1/2}\,\mathrm{d}y\\
&&{}+s^{1/2-\lambda
}\int_{s}^{t}(y-s)^{\lambda-3/2}y^{\lambda-1/2}\bigl(\log(y-s)+\log y\bigr)\,\mathrm{d}y.
\end{eqnarray*}
Then, for every $T>0$, there exists a constant $C_{a,b,T}$ such that
for every $ \lambda\in[a,b]$ and $s\in(0,T]$,
\[
\biggl| \frac{
\partial}{\partial\lambda}K_{\lambda}(t,s)\biggr| \leq
C_{a,b,T}(1\vee|
\log s| )s^{ {1/2} - b}=: \Phi_T(s),
\]
which concludes the proof.
\end{pf}

The following proposition gives the covariance of this
process.
\begin{proposition}
\label{propcov}
Let $X=\{X(t,\lambda),t\geq0,\lambda\in( {1/2}
,1)\} $ be the two-parameter process given by
$X(t,\lambda)=\int_{0}^{t}K_{\lambda}(t,u)W(\mathrm{d}u)$. Then
\[
\mathbb{E}[X(t,\lambda)X(s,\lambda^{\prime
})]=\int_{0}^{t}\mathrm{d}y\int_{0}^{s}\mathrm{d}z\,\widetilde{\beta}(y,z,\lambda,\lambda
^{\prime})| y-z| ^{\lambda+\lambda^{\prime}-2}
\biggl( \frac{y}{z}\biggr)^{\lambda-\lambda^{\prime}} ,
\]
where
\[
\widetilde{\beta}(y,z,\lambda,\lambda^{\prime})=\beta(2-\lambda
-\lambda^{\prime},\lambda^{\prime}- {1/2} )1_{\{y>z\}}+\beta
(2-\lambda-\lambda^{\prime},\lambda- {1/2} )1_{\{y<z\}}
\]
and $\beta(a,b)$ $(a,b>0)$ is the beta function. In particular,
\[
\mathbb{E}[X(t,\lambda)^{2}]  =  \int_{0}^{t}\mathrm{d}y\int_{0}^{t}\mathrm{d}z\,\beta(2-2\lambda,\lambda- {1/2} )| y-z| ^{2\lambda-2}
 =
{t^{2\lambda}\over c_{\lambda}^2}
\]
with
%
\begin{equation}
\label{defcH}
c_{\lambda}^2 =\frac{2\curpi\lambda(\lambda-1/2)^3}{\Gamma(2-2\lambda
)\Gamma(\lambda+{1/2}
)^{2}\sin(\curpi(\lambda-{1/2}))} .
\end{equation}
\end{proposition}
\begin{pf}
We have
\begin{eqnarray*}
&&\mathbb{E}[X(t,\lambda)X(s,\lambda^{\prime
})] \\
&&\quad = \int_{0}^{t\wedge s}K_{\lambda}(t,u)K_{\lambda^{\prime
}}(s,u)\,\mathrm{d}u\\
&&\quad =  \int_{0}^{t\wedge s}u^{1-\lambda-\lambda^{\prime}}\biggl( \int
_{u}^{t}(y-u)^{\lambda- {3/2} }y^{\lambda- {1/2}
}\,\mathrm{d}y\biggr)\biggl( \int_{u}^{s}(z-u)^{\lambda^{\prime}-
{3/2} }z^{\lambda^{\prime}-
{1/2} }\,\mathrm{d}z\biggr)\,\mathrm{d}u\\
&&\quad = \int_{0}^{t}\mathrm{d}y\int_{0}^{s}\,\mathrm{d}z  y^{\lambda- {1/2} }z^{\lambda
^{\prime
}- {1/2} }\int_{0}^{y\wedge z}u^{1-\lambda-\lambda
^{\prime
}}(y-u)^{\lambda- {3/2} }(z-u)^{\lambda^{\prime}- {3/2} }\,\mathrm{d}u.
\end{eqnarray*}
We fix $y>z$ and calculate the following integral by making the successive
substitutions $u=vz$, $w={(y-vz)/(1-v)},$ $t={w/y}$ and $s={1/t}$:
\begin{eqnarray*}
&&\int_{0}^{z}u^{1-\lambda-\lambda^{\prime}}(y-u)^{\lambda- {3/2}
}(z-u)^{\lambda^{\prime}- {3/2} }\,\mathrm{d}u\\
&&\quad=z^{ {1/2} -\lambda}\int_{0}^{1}\lambda^{\prime
}(y-vz)^{\lambda-
{3/2} }(1-v)^{\lambda^{\prime}- {3/2} }\,\mathrm{d}v\\
&&\quad=z^{ {1/2} -\lambda}(y-z)^{\lambda+\lambda^{\prime
}-2}\int_{y}^{+\infty}(w-y)^{1-\lambda-\lambda^{\prime}}w^{\lambda-
{3/2} }\,\mathrm{d}w\\
&&\quad=z^{ {1/2} -\lambda}y^{ {1/2} -\lambda^{\prime
}}(y-z)^{\lambda
+\lambda^{\prime}-2}\int_{1}^{+\infty}(t-1)^{1-\lambda-\lambda
^{\prime
}}t^{\lambda- {3/2} }\,\mathrm{d}t\\
&&\quad=z^{ {1/2} -\lambda}y^{ {1/2} -\lambda^{\prime
}}(y-z)^{\lambda
+\lambda^{\prime}-2}\beta(2-\lambda-\lambda^{\prime},\lambda
^{\prime
}- {1/2} ).
\end{eqnarray*}
In the same way, for $y<z$, we obtain
\begin{eqnarray*}
&&\int_{0}^{y}u^{1-\lambda-\lambda^{\prime}}(y-u)^{\lambda- {3/2}
}(z-u)^{\lambda^{\prime}- {3/2} }\,\mathrm{d}u\\
&&\quad=z^{ {1/2} -\lambda}y^{{1}/{2}-\lambda^{\prime}}(z-y)^{\lambda+\lambda^{\prime}-2}\beta
(2-\lambda-\lambda^{\prime},\lambda- {1/2} ).
\end{eqnarray*}
This concludes the proof.
\end{pf}

In the sequel, we need the estimates which we establish in the
following proposition.
\begin{proposition}
\label{increments}
The process $X$ satisfies the following estimates:
\begin{longlist}[(a)]
\item[(a)] for all $s$ and $t\geq0$
%
\begin{equation}
\label{increments1}
\mathbb{E}\bigl[\bigl(X(t,\lambda)-X(s,\lambda
)\bigr)^{2}\bigr]=c_{\lambda}^{-2}| t-s|^{2\lambda} ;
\end{equation}
\item[(b)] for every $T>0$, there exists a constant $C_T>0$ such that for every
$t\in[0,T]$ and every $\lambda$ and $\lambda^{\prime}\in[a,b]$,
%
\begin{equation}
\label{increments2}
\mathbb{E}\bigl[\bigl(X(t,\lambda)-X(t,\lambda^{\prime})\bigr)^{2}\bigr]\leq
C_{T}|\lambda-\lambda^{\prime}|^{2} .
\end{equation}
\end{longlist}
\end{proposition}
\begin{pf}(a) For every $\lambda\in[a,b]$, the process $X(\cdot, \lambda)\dvtx t\mapsto
X(t , \lambda)$ is a fractional Brownian motion with variance
$c_{\lambda}^{-2}$, so we deduce (\ref{increments1}).

(b) We have
\[
\mathbb{E}\bigl[\bigl(X(t,\lambda)-X(t,\lambda^{\prime
})\bigr)^{2}\bigr]=\int_{0}^{t}\bigl(K_{\lambda}(t,u)-K_{\lambda^{\prime}}(t,u)\bigr)^{2}\,\mathrm{d}u.
\]
There exists $\xi=\xi(\lambda,\lambda^{\prime}) \in[\min\{
\lambda,\lambda^{\prime}\},\max
\{\lambda,\lambda^{\prime}\}]$ such that
\[
K_{\lambda}(t,u)-K_{\lambda^{\prime}}(t,u)=( \lambda-\lambda
^{\prime})
\biggl| \frac{\partial}{\partial
\lambda}K_{\lambda}(t,u)\biggr| _{\lambda=\xi=\xi(\lambda
,\lambda^{\prime
})}.
\]
Then, thanks to Lemma \ref{estimateKH}, for every $t$, $\lambda$ and
$\lambda^{\prime}$, we obtain
\[
\mathbb{E}\bigl[\bigl(X(t,\lambda)-X(t,\lambda^{\prime
})\bigr)^{2}\bigr]\leq|\lambda- \lambda^{\prime}|^2\int_{0}^{T}(\Phi
_T(u))^{2}\,\mathrm{d}u ,
\]
which concludes the proof.
\end{pf}

From the last proposition, we can deduce the continuity of $B_h$.

\begin{corollary}
\label{continuity}
The process $B_h$ defined above has continuous
trajectories.
\end{corollary}
\begin{pf}
We deduce from Proposition \ref{increments} and the assumption $\sup h<
\beta$ that for every $s$ and $t$ in a compact interval $[0,T]$ such
that $|t-s| < 1$, we have
\begin{eqnarray*}
\e\bigl[\bigl(B_h(t)-B_h(s)\bigr)^2\bigr] & = & \e\bigl[\bigl(X(t,h(t))-X(s,h(s))\bigr)^2\bigr]\\
& \leq& 2\e\bigl[\bigl(X(t,h(t))-X(t,h(s))\bigr)^2\bigr]\\
& &{} +2\e\bigl[\bigl(X(t,h(s))-X(s,h(s))\bigr)^2\bigr]\\
& \leq& 2C_T|h(t)-h(s)|^2 +2c_{h(s)}^{-2}|t-s|^{2h(s)}\\
& \leq& 2C_T|t-s|^{2\beta} +2\sup_{\lambda\in[a,b]} (c_{\lambda
}^{-2})|t-s|^{2a}.
\end{eqnarray*}
Then $\e[(B_h(t)-B_h(s))^2]/|t-s|^{2a}$ is bounded and since $B_h$ is
Gaussian, we can deduce from~\cite{Bil} its continuity.
\end{pf}

We now deal with the local self-similarity property of $B_h$.

\begin{proposition}
\label{selfsimilar}
The process $B_h$ is locally self-similar. More precisely, for every
$t$, we have the following
convergence in distribution:
\[
\lim_{\varepsilon\to0} \biggl( { B_h(t+\ve u)-B_h(t) \over\ve^{h(t)}
} \biggr)_{u\geq0}
= \bigl( c_{h(t)}^{-1}B_{h(t)}(u)\bigr)_{u\geq0} ,
\]
where $\lim_{\varepsilon\to0}$ stands for the limit in distribution
in the space of continuous
functions endowed with the uniform norm on every compact set.
\end{proposition}
\begin{pf}
Let us start by proving the convergence of the finite-dimensional
distribution. Because~$B_h$ is Gaussian,
it suffices to show the convergence of the second order moments.
We then can write, for every $u$ and $v$,
\begin{eqnarray*}
&&{1\over\ve^{2h(t)}}\e\bigl[ \bigl( B_h(t+\ve u)-B_h(t)\bigr)
\bigl( B_h(t+\ve v)-B_h(t)\bigr) \bigr] \\
&&\quad = {1\over\ve^{2h(t)}}\bigl(I_1(\ve)+I_2(\ve)+I_3(\ve)+I_4(\ve
)\bigr) ,
\end{eqnarray*}
where
\begin{eqnarray*}
I_1(\ve) &=& \e\bigl[ \bigl( X\bigl(t+\ve u, h( t )\bigr)-X(t, h( t ))\bigr) \bigl( X\bigl(t+\ve v, h( t )\bigr)-X(t, h( t ))\bigr)
\bigr] , \\
I_2(\ve) &=& \e\bigl[ \bigl( X\bigl(t+\ve u, h( t+\ve u )\bigr)-X\bigl(t+\ve u , h( t
)\bigr)\bigr) \bigl( X\bigl(t+\ve v, h( t )\bigr)-X( t , h( t ))
\bigr) \bigr], \\
I_3(\ve) &=& \e\bigl[ \bigl( X\bigl(t+\ve u, h( t )\bigr)-X( t , h( t ))\bigr) \bigl( X\bigl(t+\ve v, h( t+\ve v )\bigr)-X\bigl( t+\ve v, h(
t )\bigr)\bigr) \bigr] , \\
I_4(\ve) &=& \e\bigl[ \bigl( X\bigl(t+\ve u, h( t+\ve u )\bigr)\hspace*{-0.4pt}-\hspace*{-0.4pt}X\bigl(t+\ve u , h( t
)\bigr)\bigr) \bigl( X\bigl(t+\ve v, h( t+\ve v )\bigr)\hspace*{-0.4pt}-\hspace*{-0.4pt}X\bigl( t+\ve v, h(
t )\bigr)\bigr) \bigr] .
\end{eqnarray*}
By the self-similarity of the fractional Brownian motion and the
stationarity of its increments, we obtain
%
\begin{equation}
\label{conv1}
\lim_{\ve\to0}{1\over\ve^{2h(t)}}I_1(\ve)={1\over
2c_{h(t)}^{2}}\bigl(|u|^{2h(t)}+|v|^{2h(t)}-|u-v|^{2h(t)}\bigr).
\end{equation}
Then, because of the Cauchy--Schwarz inequality, it is enough to prove
%
\begin{equation}
\label{conv2}
\lim_{\ve\to0}{1\over\ve^{2h(t)}}\e\bigl[ \bigl( X\bigl(t+\ve u, h( t+\ve
u )\bigr)-X\bigl(t+\ve u , h( t )\bigr)\bigr)^2 \bigr] =0
\end{equation}
to get
%
\begin{equation}
\label{conv3}
\lim_{\ve\to0}{1\over\ve^{2h(t)}}\bigl(I_2(\ve) + I_3(\ve) + I_4(\ve) \bigr)=0.
\end{equation}
Thus, using Lemma \ref{estimateKH}, we obtain
\begin{eqnarray*}
&& {1\over\ve^{2h(t)}}\e\bigl[ \bigl( X\bigl(t+\ve u, h( t+\ve u )\bigr)-X\bigl(t+\ve
u , h( t )\bigr)\bigr)^2 \bigr] \\
&&\quad = {1\over\ve^{2h(t)}} \int_{0}^{t+\ve u}
\bigl( K_{h( t+\ve u )}(t+\ve u, s)-K_{h( t )}(t+\ve u, s)\bigr)^2\,\mathrm{d}s\\
&&\quad \leq { (h(t+\ve u) - h(t) )^2 \over\ve^{2h(t)}} \int
_{0}^{t+\ve u}
\sup_H\biggl| {\partial\over\partial H}K_{H}(t+\ve u, s)\biggr|^2\,\mathrm{d}s\\
&&\quad \leq { C\ve^{2\beta-2h(t)} } \int_{0}^{t+u}
| \Phi_{t+u}(s)|^2 \,\mathrm{d}s
\mathop{{\longrightarrow}}_{\ve\to0} 0.
\end{eqnarray*}
It now remains to prove tightness in the space of continuous functions
endowed with the uniform norm.
We also consider $T>0$ such that $t$, $t+\ve u$ and $t+\ve v\in[0,T]$
for all $\varepsilon$. Performing similar calculations as in the proof
of Corollary \ref{continuity}, we get that there exist $C_T>0$ such that
\begin{eqnarray*}
&& \e\biggl[\biggl( { B_h(t+\ve u)-B_h(t) \over\ve^{h(t)} }- { B_h(t+\ve
v)-B_h(t) \over\ve^{h(t)} } \biggr)^2\biggr]\\
&&\quad = {1\over\ve^{2h(t)}}\e\bigl[\bigl( { B_h(t+\ve u)} - {
B_h(t+\ve v) } \bigr)^2\bigr]\\
&&\quad \leq {C_T\over\ve^{2h(t)}}| {\ve u - \ve v }
|^{2h(t+\ve v)}\\
&&\quad = {C_T \ve^{2(h(t+\ve v)-h(t))}}| { u - v }
|^{2h(t+\ve v)} .
\end{eqnarray*}
Since $h$ is $\beta$-H\"older, it follows that $\ve^{2(h(t+\ve
v)-h(t))}$ is uniformly bounded. Moreover, $h(t+\break \ve v)\geq a$, thus
\[
 \e\biggl[\biggl( { B_h(t+\ve u)-B_h(t) \over\ve^{h(t)} }- { B_h(t+\ve
v)-B_h(t) \over\ve^{h(t)} } \biggr)^2\biggr] \leq {C_T }| { u -
v } |^{2a}.
\]
This completes the proof.
\end{pf}

It is classical to deduce pointwise H\"older continuity from local
self-similarity \cite{BCI}. We recall that the pointwise H\"older
continuity of a function $f$ is characterized by the pointwise H\"older exponent
$\alpha_{f}(t_0)$ defined at each point $t_0$ as
\[
\alpha_{f}(t_0)=\sup\biggl\{ \alpha>0 \dvtx \lim_{t\to
t_0}{|f(t)-f(t_0)|\over|t-t_0|^{\alpha}}=0 \biggr\}.
\]
\begin{proposition}
For every $t_0\in\re_+$, the pointwise H\"older exponent $\alpha_{B_h}(t_0)$
of $B_h$ is almost surely equal to $h(t_0)$.
\end{proposition}
\begin{pf}
We deduce from
Proposition \ref{selfsimilar} and \cite{BCI} that $\alpha
_{B_h}(t_0)\leq h(t_0)$.
We now prove that $\alpha_{B_h}(t_0)\geq h(t_0)$. By Proposition \ref
{increments},
there exists $C>0$ such that for every $\ve>0$ and every
$s$ and $t\in[t_0-\ve, t_0+\ve]$ such that $|t-s|<1$, we have
\[
\e\bigl[\bigl(B_h(t)-B_h(s)\bigr)^2\bigr] \leq C |t-s|^{2\inf_{[ t_0-\ve, t_0+\ve]}h}.
\]
By using the fact that $B_h$ is Gaussian and applying Kolmogorov's
theorem \cite{Bil},
we get that $\lim_{t\to t_0}{|B_h(t)-B_h(t_0)|/ |t-t_0|^{\alpha}}=0 $
for every $\alpha< \inf_{[ t_0-\ve, t_0+\ve]}h$.
This holds for every $\ve>0$, so by continuity of $h$, we have $\lim
_{t\to t_0}{|B_h(t)-B_h(t_0)|/ |t-t_0|^{\alpha}}=0 $ for every $\alpha
< h(t_0)$. We can deduce that
$\alpha_{B_h}(t_0)\geq h(t_0)$ and hence $\alpha_{B_h}(t_0)=h(t_0)$.
\end{pf}

\section{Stochastic calculus for $B_{h}$}

The aim of this section is to apply the stochastic calculus developed
by Al\`{o}s,
Mazet and Nualart in \cite{AMN} to get a stochastic
integral for $B_{h}$ and an It\^{o} formula. We recall that in \cite
{AMN}, the following
hypothesis, called (K4),
appears for regular kernels:
\begin{itemize}
\item[$\bullet$] \textbf{Hypothesis} (K4). For all $s\in\lbrack0,T)$, $K(\cdot
,s)$ has bounded variation on the interval $(s,T]$ and $\int_{0}^{T}| K|
((s,T],s)^{2}\,\mathrm{d}s<\infty$ (where, for all $s\in\lbrack0,T)$, $
|K|((s,T],s)$ denotes the total variation of $K(\cdot,s)$ on $(s,T]$).
\end{itemize}
\begin{lemma}
Suppose that $h$ is of bounded variation on $(s,T]$ for
all $s\in\lbrack0,T).$ Then \textup{(K4)} holds for $(t,s)\mapsto
K_{h(t)}(t,s)$ defined by
(\ref{defKh}).
\end{lemma}
\begin{pf}
Let
\[
\operatorname{Var}_{(s,T]}^{n}(\cdot,s)=\sup_{
t_{0}=s<t_{1}<\cdots<t_{n}=T}\sum_{i=1}^{n}\bigl|
K_{h(t_{i})}(t_{i},s)-K_{h(t_{i-1})}(t_{i-1},s)\bigr|
\]
and suppose, without
restriction of generality, that $T=1$.
Then
\[
\bigl| K_{h(t_{i})}(t_{i},s)-K_{h(t_{i-1})}(t_{i-1},s )\bigr| \leq
I_{1}(i)+I_{2}(i),
\]
where
\[
I_{1}(i) = \bigl|K_{h(t_{i})}(t_{i},s)-K_{h(t_{i})}(t_{i-1},s)\bigr|
\]
and
\[
I_{2}(i) = \bigl|K_{h(t_{i})}(t_{i-1},s)-K_{h(t_{i-1})}(t_{i-1},s) \bigr| .
\]
We have
\begin{eqnarray*}
I_{1} (i)& = & s^{ {1/2} -h(t_{i})}\int
_{t_{i-1}}^{t_{i}}(y-s)^{h(t_{i})- {3/2} }y^{h(t_{i})- {1/2} }\,\mathrm{d}y\\
&\leq& s^{ {1/2} -b}\int_{t_{i-1}}^{t_{i}}(y-s)^{a- {3/2} }y^{a-
{1/2} }\,\mathrm{d}y.
\end{eqnarray*}
Therefore,
%
\begin{eqnarray}
\label{estimateI1}
&&\sum_{i=1}^{n}\bigl|
K_{h(t_{i})}(t_{i},s)-K_{h(t_{i})}(t_{i-1},s)\bigr|\nonumber\\
&&\quad  \leq
s^{ {1/2} -b}\int_{s}^{1}(y-s)^{a- {3/2} }y^{a-{1/2}}\,\mathrm{d}y
\\
&&\quad\leq C(a)s^{ {1/2} -b}.\nonumber
\end{eqnarray}
Regarding $I_{2}$, we have
\begin{eqnarray*}
&&\sum_{i=1}^{n}\bigl|
K_{h(t_{i})}(t_{i-1},s)-K_{h(t_{i-1})}(t_{i-1},s)\bigr|\\
&&\quad
\leq\sum_{i=1}^{n}| h(t_{i})-h(t_{i-1})| \sup_{s\leq t
\leq1,a\leq\lambda\leq b}\biggl|\frac{\partial}{\partial
\lambda}K_{\lambda
}(t,s)\biggr| \\
&&\quad
\leq\mathop{\Var}_{[0,T]}h \sup_{s\leq t
\leq1,a\leq\lambda\leq b}\biggl|\frac{\partial}{\partial
\lambda}K_{\lambda
}(t,s)\biggr| .
\end{eqnarray*}
The proof of Lemma\vspace*{1pt} \ref{estimateKH} implies that $\operatorname{Var}_{(s,T]}^{n}(\cdot,s)\leq C(1\vee|\log s|)s^{ {1/2} -b},$ where the constant $C>0$ depends on $h$ (but
not on $n$)$.$ This implies that (K4) holds for $K_{h}$.
\end{pf}
\begin{remark}
Since $h$ is also supposed $\beta$-H\"{o}lder continuous
for some $\beta\leq1,$ we will hereafter suppose that $h$ is
Lipschitz continuous. Note that $\lim_{t\searrow u}%
K_{h(t)}(t,u)=0$ for all $u>0.$
\end{remark}

For simplicity, we write $K$ instead of $K_{h}$, but keep in mind that $
(t,s)\rightarrow K(t,s)$ means $(t,s)\rightarrow K_{h(t)}(t,s)$ and that
differentials with respect to $t$ also act via $h$.

In the sequel, we need the following proposition regarding the variance
of $B_h$.
\begin{proposition}
The variance $s\mapsto R_{s}=\mathbb{E}[B_{h}(s)^{2}]$
is of bounded variation on $(0,T].$
\end{proposition}
\begin{pf}
We have
\begin{eqnarray*}
\sum_{i=1}^{n} | R_{s_{i+1}}-R_{s_{i}}| & = & \sum
_{i=1}^{n}\biggl|
\int_{0}^{s_{i+1}}K(s_{i+1},r)^{2}\,\mathrm{d}r-\int
_{0}^{s_{i}}K(s_{i},r)^{2}\,\mathrm{d}r\biggr| \\
&\leq
&\sum_{i=1}^{n}\int_{0}^{s_{i+1}}K(s_{i+1},r)^{2}\,\mathrm{d}r+\sum_{i=1}^{n}\biggl|
\int_{0}^{s_{i}}[K(s_{i+1},r)^{2}-K(s_{i},r)^{2}]\,\mathrm{d}r\biggr|.
\end{eqnarray*}
The functions $| K(s,r)1_{[0,s]}(r)| $ are bounded by the
square-integrable function $k(r)=| K(T,r)| + |K| ((r,T],r).$ Hence, the
first term above is upper bounded by $\int_{0}^{T}k(r)^{2}\,\mathrm{d}r.$ The second
term is upper bounded by
\begin{eqnarray*}
&&\sum_{i=1}^{n}\int_{0}^{s_{i}} |K| ((s_{i},s_{i+1}],r)|
K(s_{i+1},r)+K(s_{i},r)| \,\mathrm{d}r \\
&&\quad\leq2\int_{0}^{T}k(r)| K |((r,T],r)\,\mathrm{d}r<\infty.
\end{eqnarray*}
This concludes the proof.
\end{pf}

For any $f\in L^{2}([0,T])$, let $Kf$ be defined by ($Kf)(t)=
\int_{0}^{t}K(t,s)f(s)\,\mathrm{d}s$. Let $\mathcal{E}$ be the set of step
functions on
$[0,T]$ and let the operator $K^{\ast}$ be defined on $\mathcal{E}$ by $
(K^{\ast}\varphi)(s)=\int_{s}^{T}\varphi(t)K(\mathrm{d}t,s).$ Then $K^{\ast
}$ is
the adjoint of $K.$ In fact, for $a_{i}\in
\mathbb{R}
$, $0=s_{1}<s_{2}<\cdots<s_{n+1}=T$, $\varphi
=\sum_{i=1}^{n}a_{i}1_{(s_{i},s_{i+1}]}(s)$ and $f\in L^{2}([0,T]),$ we write
\[
(K^{\ast}\varphi
)(s)=\sum_{i=1}^{n}1_{(s_{i},s_{i+1}]}(s)
\sum_{j=i}^{n}a_{j}\bigl(K(s_{j+1},s)-K(s_{j},s)\bigr)
\]
and
%
\begin{eqnarray}
\label{5}
\int_{0}^{T}(K^{\ast}\varphi)(s)f(s)\,\mathrm{d}s
&=&\sum_{j=1}^{n}a_{j}\int_{0}^{T}
\sum_{i=1}^{j}1_{(s_{i},s_{i+1}]}(s)\bigl(K(s_{j+1},s)-K(s_{j},s)\bigr)f(s)\,\mathrm{d}s
\nonumber
\\[-8pt]
\\[-8pt]
\nonumber
&=&\sum_{j=1}^{n}a_{j}[(Kf)(s_{j+1})-(Kf)(s_{j})]=\int_{0}^{T}\varphi
(t)(Kf)(\mathrm{d}t).
\end{eqnarray}
As usual, the reproducing kernel Hilbert space (RKHS) $\mathcal{H}$ is
defined as the closure of the linear span of the indicator functions \{1$
_{[0,t]},t\in\lbrack0,T]$\} with respect to the scalar product
$\langle
1_{[0,t]},1_{[0,s]}\rangle_{\mathcal{H}}=\mathbb{E}[B_{h(t)}(t)B_{h(s)}(s)]
\equiv R(t,s).$ By replacing $f(s)\,\mathrm{d}s$ in (\ref{5}) by $W(\mathrm{d}s),$ we have,
by (1),
\[
B_{h}(\varphi)\equiv\int_{0}^{T}\varphi(t)B_{h}(\mathrm{d}t)=\int
_{0}^{T}(K^{\ast
}\varphi)(s)W(\mathrm{d}s).
\]
Therefore,
\[
\| \varphi\| _{\mathcal{H}}^{2}=\mathbb{E}[B_{h}(\varphi
)^{2}]=\| K^{\ast}\varphi\| _{L^{2}([0,T])}^{2}\leq
\int_{0}^{T}\biggl[ \int_{0}^{T}| \varphi(t)| | K|
(\mathrm{d}t,s)\biggr]^{2}\,\mathrm{d}s=:\| \varphi\| _{K}^{2}.
\]

Let us denote by $\mathcal{H}_{K}$ the completion of $\mathcal{E}$
with respect to the $\| \cdot\| _{K}$-norm. Then $\mathcal{H}
_{K}$ is continuously embedded in $\mathcal{H}.$

In order to define the stochastic integral with respect to $B_{h}$, let us
denote by $\mathcal{S}$ the set of smooth cylindrical random variables of
the form $F=f(B_{h}(\varphi_{1}),B_{h}(\varphi_{2}),\ldots
,B_{h}(\varphi
_{n})),$ where $n\geq1,$ $f\in C_{b}^{\infty}(
\mathbb{R}^{n})$ ($f$ and all its derivatives are bounded) and $\varphi
_{1},\varphi
_{2},\ldots,\varphi_{n}\in\mathcal{H}.$ Let us also denote by $
\mathbb{D}^{1,2}(\mathcal{H}_{K})$ the closure of \{$F\in\mathcal
{S}\dvtx F\in
L^{2}(\Omega,\mathcal{H}_{K}),DF\in L^{2}(\Omega\times\mathcal{H}_{K},
\mathcal{H}_{K})$\}. Then $\mathbb{D}^{1,2}(\mathcal{H}_{K})$ is
included in the domain $\operatorname{Dom}(\delta^{B_{h}})$ of the
divergence operator of $%
B_{h},$ and the integral of $u\in\operatorname{Dom}(\delta^{B_{h}})$ with
respect to
$B_{h}$ is given by
\[
\delta^{B_{h}}(u)\equiv\int_{0}^{T}u  \,\delta B_{h}=\int
_{0}^{T}(K^{\ast
}u)  \,\delta W\equiv\int_{0}^{T}\biggl[ \int_{s}^{T}u(r){K}(\mathrm{d}r,s)
\biggr] \,\delta W(s),
\]
where the last two integrals are the divergence integrals with respect to
Brownian motion. Let us recall that the integral $\delta^{B_{h}}(u)$ is
defined, for any $u\in L^{2}(\Omega,\mathcal{H}),$ as the unique element
in $L^{2}(\Omega)$ which satisfies the duality relationship $\mathbb{E}
(\delta^{B_{h}}(u)F)=\mathbb{E}\langle DF,u\rangle_{\mathcal{H}}$ for all
$F\in\mathcal{S}.$

The following It\^{o} formula, due to Al\`{o}s, Mazet and Nualart \cite
{AMN}, will
be applied in the next section. Let $F$ be a function belonging to the class
$C^{2}(\mathbb{R})$ and satisfying the condition
\[
\max\{| F(x)| ,| F^{\prime}(x)| ,| F^{\prime\prime}(x)| \}\leq
c\mathrm{e}^{\lambda| x| ^{2}},
\]
where $c$ and $\lambda$ are positive constants such that $\lambda
<\frac{1}{
4}( {\sup_{0\leqq t\leq T}}R_{t})^{-1}.$ This implies that the
process $F^{\prime}(B_{h})$ belongs to $\mathbb{D}^{1,2}(
\mathcal{H}_{K})$. The integral $\int_{0}^{t}F^{\prime}(B_{h}(s))\,\delta
B_{h}(s)$ is therefore well defined for all $t\in\lbrack0,T]$ and the
following It\^{o}-type formula holds (\cite{AMN}, Theorem 2):
%
\begin{equation}
\label{ito}
F(B_{h}(t))=F(0)+\int_{0}^{t}F^{\prime}(B_{h}(s))\,\delta B_{h}(s)+
{1\over2}
\int_{0}^{t}F^{\prime\prime}(B_{h}(s))\,\mathrm{d}R(s).
\end{equation}

\section{Local time and Tanaka formula for $B_{h}$}

First, we prove, by means of a criterion for Gaussian processes due to
Berman, that $B_{h}$ has a local time with respect to the Lebesgue measure.
We then derive a Tanaka-type formula from the It\^{o} formula of
Section 3
and show that $B_{h}$ satisfies the LND property. This implies continuity
and H\"older regularities in space and in time of the
trajectories of the local time.

\begin{definition}
For any Borel set
$C\subset\mathbb{R}_{+}$, \textit{the occupation measure} $m_{C}$ of
$B_{h}$on $C$
is defined, for
all Borel sets $A\subset\mathbb{R}$, by
$m_{C}(A)=\lambda\{t\in C,$ $B_{h(t)}(t)\in A\},$ where $\lambda$ is
the Lebesgue measure on
$\mathbb{R}_{+}$. If $m_{C}$ is absolutely continuous with respect to
the Lebesgue
measure on $
\mathbb{R}$, the \textit{local time} (or occupation density) of $B_h$ on
$C$ is
defined as the Radon--Nikodym derivative of $m_{C}$ and will be denoted by
$\{L(C,x),x\in\mathbb{R}\}$. Sometimes we write $L(t,x)$ instead of
$L([0,t],x).$
\end{definition}

This definition implies that the local time of $B_{h}$ satisfies the
occupation density formula
%
\begin{eqnarray}
\label{localtime}
\int_{C}g(t,B_{h}(t))\,\mathrm{d}t=\int_{C\times\mathbb{R}}g(t,x)L(\mathrm{d}t,x)\,\mathrm{d}x
\end{eqnarray}
for all continuous functions with compact support
$g\dvtx C\times\mathbb{R}\rightarrow\mathbb{R}$.
If $g$ does not depend explicitly on $t$, then we get the more classical
occupation density formula where the right-hand side of (\ref{localtime})
is replaced by $\int_{\re}g(x)L(C,x)\,\mathrm{d}x$.
\begin{proposition}
The local time of $B_{h}$ exists $P$-a.s.~on any
interval $[0,T]$ and is a square-integrable function of $x$.
\end{proposition}
\begin{pf}
By \cite{Ber69}, for any continuous and zero-mean Gaussian
process $\{X_{t},$ $t\in[0,T]\}$ with bounded covariance function, the
condition
%
\begin{equation}
\label{berman}
\int_{0}^{T}\!\!\!\int_{0}^{T}{\mathrm{d}s  \,\mathrm{d}t\over\sqrt{\mathbb
{E}[|X_{t}-X_{s}|^{2}]}} <+\infty
\end{equation}
is sufficient for the local time of $X$ to exist and to be a
square-integrable function of $x$.
For $(s,t)$ away from the diagonal, we write, for $s<t$,
\[
\mathbb{E}[|B_{h}(t)-B_{h}(s)|^{2}]=\int
_{0}^{s}\bigl(K_{h(t)}(t,u)-K_{h(s)}(s,u)\bigr)^{2}\,\mathrm{d}u+\int_{s}^{t}K_{h(t)}(t,u)^{2}\,\mathrm{d}u
\]
and deduce from (1) that the second term stays strictly positive as $(s,t)$
varies in $[0,t-\varepsilon]\times [0,T].$ If $(s,t)$ is close to
the diagonal, say $0\leq t-s<\varepsilon,$ then, by Proposition 5,
\[
\mathbb{E}[|B_{h}(t)-B_{h}(s)|^{2}]\sim c_{h(t)}^{-2}(t-s)^{2h(t)}
\]
as $s\nearrow t$
and a direct calculation shows that (\ref{berman}) is satisfied.
\end{pf}

Let us now derive a Tanaka-type formula for $B_{h}.$ Since the last term
in the It\^o formula (\ref{ito}) is an integral with
respect to
the variance function $R$ and since $R$ is not, in general,
increasing, but only of finite variation, the trajectorial representation
of the local time is more delicate than for Brownian motion
or for fractional Brownian motion. In fact, on the time intervals where
$R$ is
(strictly) increasing or decreasing, this formula gives an occupation
density $\widehat{L}$ related to $L$.

\begin{theorem}\label{th12}
Suppose that $h$ is continuously differentiable. Then,
for all $a\in\mathbb{R}$,
\[
|B_{h}(t)-a| -| B_{h}(s)-a|
=\int_{s}^{t}\sign\bigl(B_{h}(u)-a\bigr)\,\delta B_{h}(u)+\widehat{L}([s,t],a)
\]
$P$-a.s., where $\widehat{L}([s,t],a)=\int_{s}^{t}R^{\prime
}(u)L(\mathrm{d}u,a)$ and $ L$ is
the local time with respect to the Lebesgue measure of $B_{h}.$ On the time
intervals $[s,t]$ where $R$ is strictly increasing (resp., strictly
decreasing), $\widehat{L}$ (resp., $-\widehat{L}$) is the (positive)
occupation density of $B_{h}$ with respect to the measure induced by $R.$
\end{theorem}
\begin{remark}
\textup{(a)} No information on the local time of $B_{h}$ can be
deduced from the above formula if $R^{\prime}(u)=0$ on the interval $[s,t]$.
In fact,
\[
R^{\prime}(u)=\frac{\mathrm{d}}{\mathrm{d}u}\mathbb{E}\bigl[B_{h(u)}^{2}(u)\bigr]=\frac{\mathrm{d}}{\mathrm{d}u}
\bigl(u^{2h(u)}\bigr)=2\biggl( h^{\prime}(u)\log u+\frac{h(u)}{u}\biggr) u^{2h(u)}
\]
and
$
R^{\prime}(u)=0
$
if $h(u)=1/\log u\in(1/2,1)$ on an interval $(s,t).$ In this case, $
\widehat{L}([s,t],a)=0$ for all $a$.

\textup{(b)} Tanaka formulas for fractional Brownian motion have been shown by several
authors. We refer to the survey \cite{C} for references.
\end{remark}
\begin{pf*}{Proof of Theorem \protect\ref{th12}} For $\varepsilon>0$, let $p_{\varepsilon
}(x)=(2\curpi\varepsilon)^{-1/2}\exp(-x^2/(2\varepsilon)).$ We apply
the It\^{o} formula of the previous section to
the function
\[
F_{\varepsilon}(x)=\int_{0}^{t}F_{\varepsilon}^{\prime}(y)\,\mathrm{d}y ,
\]
where
\[
F_{\varepsilon}^{\prime}(x)=2\int_{-\infty}^{x}p_{\varepsilon
}(y)\,\mathrm{d}y-1.
\]
$F_{\varepsilon}^{\prime}(x)$ then converges to $\sign(x)$ and
$F_{\varepsilon}(x)$ converges to $| x | $ as $\varepsilon\rightarrow
0.$ By (\ref{ito}), for
$\varepsilon>0$ fixed,
%
\begin{eqnarray}
\label{ito2}
F_{\varepsilon}\bigl(B_{h}(t)-a\bigr) & = & F_{\varepsilon
}(-a)+\int_{0}^{t}\int_{s}^{t}F_{\varepsilon}^{\prime
}\bigl(B_{h}(r)-a\bigr)K(\mathrm{d}r,s)\,\delta W_{s}
\nonumber
\\[-8pt]
\\[-8pt]
\nonumber
&&{}+\int_{0}^{t}p_{\epsilon
}\bigl(B_{h}(s)-a\bigr)\,\mathrm{d}R(s).
\end{eqnarray}
Note that, by (K4), the process $\{F_{\varepsilon}^{\prime}(B_{h}(r)-a),
r\in[ 0,t]\}$ belongs to $L^2(\Omega,\mathcal{H}_{K})$ and therefore
to Dom($\delta^{B_{h}}$).
Or, equivalently, the process $\{\int_{s}^{t}F_{\varepsilon}^{\prime
}(B_{h}(r)-a)K_{h}(\mathrm{d}r,s),s\in[ 0,t]\}$ belongs to Dom($\delta^{W}$).
Clearly, $F_{\varepsilon}(B_{h}(t)-a)$ converges to
$|B_{h}(t)-a| $ in $L^2(\Omega)$ and $F_{\epsilon}(-a)$ converges to
$| a | $ if
$\varepsilon\rightarrow0$. Moreover, the process $\{F_{\varepsilon
}^{\prime}(B_{h}(r)-a),r\in[ 0,t]\}$ converges, as $\varepsilon
\rightarrow0,$ to $\{\sign(B_{h}(r)-a),r\in[ 0,t]\}$ in $L^2(\Omega
,\mathcal{H}_{K})$.
In fact, by (K4),
\[
\mathbb{E}\biggl[ \int_{0}^{t}\biggl( \int_{s}^{t}\bigl| F_{\varepsilon
}^{\prime
}\bigl(B_{h}(r)-a\bigr)-\sign\bigl(B_{h}(r)-a\bigr)\bigr| | K_{h}| (\mathrm{d}r,s)\biggr)^2\,\mathrm{d}s\biggr]
\]
is upper bounded independently of $\varepsilon$ and we may apply the
dominated convergence theorem (\cite{CNT}, Lemma 1).
Therefore, the last term in (\ref{ito2}) converges in
$L^2(\Omega)$.
Let us denote the limit by $\Lambda_{t}^{a}.$
Therefore, for any continuous function $f$ with compact support in
$\mathbb{R}$,
%
\begin{equation}
\label{integ}
\int\biggl(\int_{0}^{t}p_{\varepsilon}\bigl(B_{h}(s)-a\bigr)\,\mathrm{d}R(s)\biggr) f(a)\,\mathrm{d}a
\end{equation}
converges in $L^{1}(\Omega)$ to $\int\Lambda_{t}^{a}f(a)\,\mathrm{d}a.$ In fact,
the dominated convergence theorem applies because
\[
\int_{0}^{t}\mathbb{E}\bigl[p_{\varepsilon
}\bigl(B_{h}(s)-a\bigr)\bigr]\,\mathrm{d}R(s)=\int_{0}^{t}p_{R(s)+\varepsilon}(a)\,\mathrm{d}R(s)\leq
\int_{0}^{t}s^{-b}R^{\prime}(s)\,\mathrm{d}s<\infty.
\]
But (\ref{integ}) also converges to
$\int_{0}^{t}f(B_{h}(s))\,\mathrm{d}R(s)=\int_{0}^{t}f(B_{h}(s))R^{\prime}(s)\,\mathrm{d}s,$
where we use the fact that~$R$ is differentiable if $h$ is differentiable$.$
Hence,
\[
\int_{0}^{t}f(B_{h}(s))R^{\prime}(s)\,\mathrm{d}s=\int\Lambda_{t}^{a}f(a)\,\mathrm{d}a.
\]
By the occupation density formula (\ref{localtime}) applied to
$g(s,x)=f(x)R^{\prime}(s)$,
we get $\Lambda_{t}^{a}=\int_{0}^{t}R^{\prime}(s)L(\mathrm{d}s,a)=\widehat{L}(t,a)$
for $\lambda$-a.e.~$a,$ $P$-a.s. We can extend this to all
$a\in\mathbb{R}$ since there exists a jointly continuous version of
$L,$ and therefore of
$\widehat{L}$, as will be shown next by Berman's methods, which are
independent of the Tanaka formula.
\end{pf*}

We now state regularity properties in time and space of the
trajectories of
$L.$ The regularity properties of $\widehat{L}$ follow easily since
$\widehat{L}(t,x)=\int_{0}^{t}R^{\prime}(s)L(\mathrm{d}s,x)$. In order to show
the existence
of a jointly continuous version of $L$, we recall the hypotheses
introduced in
\cite{BDG} and show that they are satisfied for
$B_{h}.$ We recall them in terms of any real-valued separable random
process $\{X(t),t\in[0,T]\}$ with Borel sample functions.
\begin{itemize}
\item[$\bullet$] \textbf{Hypothesis} $(\mathbb{A})$. There exist numbers $\rho
_{0}>0$ and $H\in(0,1)$ and a
positive function $\psi\in L^{1}(\mathbb{R})$ such that for all
$\lambda\in\mathbb{R}$ and $t,s\in[0,T]$, $0<| t-s| <\rho_{0}$, we have
\[
\biggl| \mathbb{E}\biggl[\exp\biggl(\mathrm{i}\lambda\frac{X(t)-X(s)}{|
t-s|^{H}}\biggr)\biggr]\biggr|
\leq
\psi(\lambda).
\]
\item[$\bullet$]\textbf{Hypothesis} $(\mathbb{A}_m)$. There exist positive
constants $\delta$ and $c$ (both
eventually depending on $m\geq2$) such that for all $t_{1},t_{2},\ldots,t_{m}$
with $0=:t_{0}<t_{1}<\cdots<t_{m}\leq T$ and $| t_{m}-t_{1}| <\delta$,
we have
\[
\Biggl| \mathbb{E}\Biggl[\exp\Biggl(\mathrm{i}\sum
_{j=1}^{m}u_{j}\bigl(X(t_{j})-X(t_{j-1})\bigr)\Biggr)
\Biggr]\Biggr| \leq\prod _{j=1}^{m}\bigl| \mathbb{E}\bigl[\exp
\bigl(\mathrm{i}cu_{j}\bigl(X(t_{j})-X(t_{j-1})\bigr)\bigr)\bigr]\bigr|
\]
for all $u_{1},u_{2},\ldots,u_{m}\in\mathbb{R}$.
\end{itemize}
Hypothesis $(\mathbb{A})$ is evidently satisfied for self-similar
processes with
stationary increments: here, $\psi(\lambda)=| \mathbb{E}[\mathrm{e}^{\mathrm{i}\lambda
X(1)}]|$. Hypothesis $(\mathbb{A})$ is also closely related to asymptotic
self-similarity and, in fact, it holds for a
fairly large class of processes (see Proposition \ref{verifH}).

If $X$ has independent increments, $(\mathbb{A}_m)$ is trivially true
for all
$m\geq2.$ When the left- and the right-hand sides of the inequality
are applied to
the characteristic function of a Gaussian process $X$, we get the condition
which is known in the literature as \textit{local
non-determinism} (LND). Loosely speaking, LND says how much dependence is
allowed for the increments of the process if the local time has
certain regularity properties. As a general rule, the trajectories of local
time become more regular if the trajectories of the process become less regular.
\begin{proposition}
For all $T>0$, the process $\{B_{h}(t),t\in[0,T]\}$ satisfies
hypotheses $(\mathbb{A})$
and $(\mathbb{A}_m)$ for all $m\geq2$ with $H=\max_{0\leq t\leq T}h(t)$.
\end{proposition}
\begin{pf} Let us start by showing $(\mathbb{A})$.
By similar calculations as in the proof of Proposition \ref
{selfsimilar}, for every $t\in[0,T]$,
\[
\biggl|\mathbb{E}\bigl[\varepsilon^{-2h(t)}\bigl(B_h(t+\varepsilon)-B_h(t)\bigr)
^{2}\bigr]-{1\over c_{h(t)}^{2}}\biggr|\leq\ve^{\beta-\sup_{[0,T]}h}\int
_{0}^{T+1}(\Phi_{T+1}(s))^2\,\mathrm{d}s.
\]
We can now conclude by applying Proposition \ref{verifH}.

We now prove $(\mathbb{A}_m)$ for all $m\geq2.$ We show that $B_{h}$
satisfies the
LND property as introduced by Berman for Gaussian processes.
For simplicity, we write $B$ instead of $B_{h}$ in this proof. Let
$t_{1}<t_{2}< \cdots<t_{m}$ and let $\mathcal{V}_{m}$ be the relative
conditioning error
given by
\[
\mathcal{V}_{m}=\frac{\Var[B(t_{m})-B(t_{m-1})|B(t_{1}), \ldots,B(t_{m-1})]}
{\Var[B(t_{m})-B(t_{m-1})]}.
\]
The Gaussian process is said to be LND if
%
\begin{eqnarray}
\label{berman2}
\mathop{\liminf_{c\searrow0^{+}}}_{0<t_{m}-t_{1}\leq c}\mathcal{V}_{m}>0.
\end{eqnarray}
This condition means that a small increment of the process is not completely
predictable on the basis of a finite number of observations from the
immediate past. For Gaussian processes, Berman \cite{Ber73} has proven that
(\ref{berman2}) implies $(\mathbb{A}_m)$. More precisely, he has shown
that if
$X$ satisfies~(\ref{berman2}) for
all $m\geq2$, then there exist constants $C_{m}>0$ and $\delta_{m}>0$ such
that, for all $u_{1},u_{2},\ldots,u_{m}$ $\in\mathbb{R}$,
\[
\Var\Biggl[ \sum_{j=1}^{m}u_{j}[B(t_{j})-B(t_{j-1})]\Biggr] \geq
C_{m}\sum_{j=1}^{m}u_{j}^{2}\Var[ B(t_{j})-B(t_{j-1})] ,
\]
where $t_{0}=0$ and $t_{1}<\cdots<t_{m}$ are different and lie in an interval
of length at most $\delta_{m}$. This implies $(\mathbb{A}_m)$.
In order to prove that $B$ verifies (\ref{berman2}), fix $m\geq2$ and
let $
0<t_{1}<t_{2}<\cdots<t_{m}<T$. By Proposition \ref{increments},
%
\begin{equation}
\label{newproof1}
\Var[B(t_{m})-B(t_{m-1})]\leq
C_{h(T)}(t_m-t_{m-1}) ^{2h(t_m)},
\end{equation}
where $C_{h(T)}$ is a constant depending on $T$.
We now study the conditional variance. Because of the definition of
$B_h$ in terms of $W$, we have the inclusion of $\sigma$-algebra
\[
\sigma\{ B(t_{1}),\ldots,B(t_{m-1})\} \subset\sigma\{ W(t), t\in[0,
t_{m-1}]\}.
\]
We then get
%
\begin{eqnarray}
\label{newproof2}
&&\Var[B(t_{m})-B(t_{m-1})|B(t_{1}),\ldots,B(t_{m-1})] \nonumber\\
&&\quad \geq
\Var\bigl[B(t_{m})-B(t_{m-1})| W(t), t\in[0, t_{m-1}] \bigr]\nonumber\\
&&\quad=  \Var\bigl[B(t_{m}) | W(t), t\in[0, t_{m-1}] \bigr]\\
&&\quad=  \Var\biggl[ \int_{t_{m-1}}^{t_m}K_{h(t_m)}(t_m,u)W(\mathrm{d}u) \biggr]\nonumber\\
&&\quad =  \int_{t_{m-1}}^{t_m}\bigl(K_{h(t_m)}(t_m,u)\bigr)^2  \,\mathrm{d}u
.\nonumber
\end{eqnarray}
We simplify the notation by letting $t:=t_m$, $s:=t_{m-1}$ and
$H:=h(t)$ and considering $\int_{s}^{t}(K_{H}(t,u))^2  \,\mathrm{d}u $.
We recall that we have the relation
%
\begin{equation}
\label{newproof3}
\int_{s}^{t}(K_{H}(t,u))^2  \,\mathrm{d}u=\e\bigl[\bigl(B_H(t)-B_H(s)\bigr)^2\bigr]-\int
_{0}^{s}\bigl(K_{H}(t,u) - K_{H}(s,u)\bigr)^2  \,\mathrm{d}u
\end{equation}
so we study $\int_{0}^{s}(K_{H}(t,u) - K_{H}(s,u))^2  \,\mathrm{d}u$. We make the
same substitutions as in the proof of Proposition \ref{propcov} and
thereby obtain
\begin{eqnarray*}
&&\int_{0}^{s}\bigl(K_{H}(t,u) - K_{H}(s,u)\bigr)^2  \,\mathrm{d}u\\
&&\quad=  \int_{0}^{s}u^{1-2H}\biggl( \int_{s}^{t}(y-u)^{H - {3/2} }y^{H- {1/2}
}\,\mathrm{d}y\biggr)^2  \,\mathrm{d}u\\
&&\quad =  \int_{s}^{t}\mathrm{d}y\int_{s}^{t}\mathrm{d}z \, y^{H - {1/2} }z^{H- {1/2} }\int_{0}^{y\wedge z}u^{1-2H}(y-u)^{H - {3/2}
}(z-u)^{H- {3/2} }\,\mathrm{d}u \\
&&\quad = \int_{s}^{t}\mathrm{d}y\int_{s}^{t}\mathrm{d}z \, |y-z|^{2H-2}\int_1^{\Psi
_1(s,y,z)}(v-1)^{1-2H}v^{H-3/2} \,\mathrm{d}v ,
\end{eqnarray*}
where
\[
\Psi_1(s,y,z) = {yz-zs\over yz-ys}1_{y>z}+ {yz-ys\over yz-zs}1_{y\leq z}.
\]
Then, using (\ref{newproof3}) and making the substitutions $(y,z)\to
(y+s,z+s)$ and $(y,z)\to((t-s)y,(t-s)z)$, we get
\[
\int_{s}^{t}(K_{H}(t,u))^2  \,\mathrm{d}u=(t-s)^{2H}\int_{0}^{1}\mathrm{d}y\int
_{0}^{1}\mathrm{d}z \, |y-z|^{2H-2}\int_{\Psi_2(s,t,y,z)}^{\infty
}(v-1)^{1-2H}v^{H-3/2} \,\mathrm{d}v ,
\]
where
\[
\Psi_2(s,t,y,z) = {(t-s)yz+ys\over(t-s)yz+zs}1_{y>z}+ {(t-s)yz+zs\over
(t-s)yz+ys}1_{y\leq z}.
\]
We can therefore deduce that
%
\begin{equation}
\label{newproof4}
\liminf_{(t-s)\to0}(t-s)^{-2H}\int_{s}^{t}(K_{H}(t,u))^2  \,\mathrm{d}u\geq
\mathcal{I} ,
\end{equation}
where
\[
\mathcal{I}=\int_{0}^{1}\mathrm{d}y\int_{0}^{1}\mathrm{d}z \, |y-z|^{2b-2}\int_{\Psi
_3(y,z)}^{\infty}(v-1)^{1-2b}v^{a-3/2} \,\mathrm{d}v > 0
\]
with
\[
\Psi_3(y,z)=\sup\biggl\{{y\over z}, {z\over y}, 2 \biggr\}.
\]
Combining (\ref{newproof1}), (\ref{newproof2}) and (\ref{newproof4}),
we then get
%
\begin{eqnarray}
\label{berman3}
\mathop{\liminf_{c\searrow0^{+}}}_{0<t_{m}-t_{1}\leq c}\mathcal
{V}_{m}\geq{\mathcal{I} \over C_{h(T)}}>0,
\end{eqnarray}
which concludes the proof.
\end{pf}

The beginning of the proof of Proposition 14 shows that the hypothesis
$(\mathbb{A})$ holds for
a much larger class of processes. In fact, we have only used the assumption
of the following proposition.
\begin{proposition}
\label{verifH}
Let $\{X_{t},t\in[0,T]\}$ be
a centered Gaussian process. Suppose that for some positive continuous functions
$f\dvtx [0,T]\rightarrow(0,1)$ and $g\dvtx [0,T]\rightarrow(0,\infty)$,
%
\begin{equation}
\label{condverifH}
\mathbb{E}\bigl[\varepsilon^{-2f(t)}\bigl(X(t+\varepsilon)-X(t)\bigr)
^{2}\bigr]{\longrightarrow}_{\varepsilon\rightarrow0}g(t)
\end{equation}
uniformly in$\ t$. Hypothesis $(\mathbb{A})$ then holds with $H=\sup f $.
\end{proposition}
\begin{pf} Let us fix $H=\sup f $. Because $X$ is Gaussian and
centered, we have
%
\begin{equation}
\label{expgauss}
\e\biggl[ \exp\biggl( \mathrm{i}\lambda{ X(t+\ve) - X(t) \over\ve^{H} } \biggr)
\biggr]=
\exp\biggl( -{\lambda^2\over2} \e\biggl[ \biggl( { X(t+\ve) - X(t) \over
\ve^{H} }\biggr)^2\biggr] \biggr).
\end{equation}
Because of (\ref{condverifH}), there exists $\ve_0$ such that for every
$\ve$
satisfying $|\ve|<\ve_0$ and for every $t$, we have
%
\begin{equation}\label{expgauss2}
\e\biggl[ \biggl( { X(t+\ve) - X(t) \over\ve^{f(t)} }\biggr)^2\biggr]
\geq{c\over2} ,
\end{equation}
where $c=\inf_{[0,T]} g $.
Besides,
$\ve^{2f(t)-2H}\geq1 $ and thus
%
\begin{equation}\label{expgauss4}
\e\biggl[ \biggl( { X(t+\ve) - X(t) \over\ve^{H} }\biggr)^2\biggr] =
\ve^{2f(t)-2H} \e\biggl[ \biggl( { X(t+\ve) - X(t) \over\ve^{f(t)}
}\biggr)^2\biggr]
\geq{c\over2}.
\end{equation}
Combining (\ref{expgauss}) and (\ref{expgauss4}), we then get, for
every $\lambda$, and $t$ and $s$ satisfying $|t-s|<\ve_0$, that
\[
\biggl|
\e\biggl[ \exp\biggl( \mathrm{i}\lambda{ X(t+\ve) - X(t) \over\ve^{H} } \biggr)
\biggr]\biggr|
\leq
\exp\biggl( -{\lambda^2c\over4} \biggr).
\]
We then choose $\psi(\lambda)=\exp( -{\lambda^2c / 4} )$ to
conclude the proof.
\end{pf}

Let us now state a regularity result for the trajectories of $L.$ The
following theorem states that $L(t,x)$ is H\"{o}lder continuous in $t$
of order
$1-H$ and H\"{o}lder continuous in $x$ of order $\frac{1-H}{2H}$, where~$H$
is the constant appearing in $(\mathbb{A})$. For the proofs, we refer
to \cite{BDG},
where it is shown that these regularities hold for any process starting at
zero and satisfying $(\mathbb{A})$ and $(\mathbb{A}_m)$ for all $m\geq
2$. Here, the function $\psi$ in $(\mathbb{A})$ satisfies
\[
\int_{\mid u\mid\geq1}| u|^{(1-H)/H}\psi(u)\,\mathrm{d}u<\infty.
\]
\begin{theorem}
$\{B_{h}(t),t\in[0,T]\}$ has a jointly continuous local time
$\{L(t,x),(t,x)\in[0,T]\times\mathbb{R}\}$.
Moreover, for any compact set $K\subset\mathbb{R}$
and any interval $I\subset[0,T]$ with length less than $\rho_{0}$
(the constant appearing in $(\mathbb{A})$):
\begin{longlist}[(ii)]
\item[(i)] if $0<\xi< (1-H)/{2H},$ then $| L(I,x)-L(I,y)| \leq\eta
|x-y|^{\xi}$ for all $x,y\in K$;

\item[(ii)] if $0<\delta<1-H,$ then $\sup_{x\in K}L(I,x)\leq\eta|I|^{\delta},$
where $\eta$ is a random variable, almost surely positive and finite, and
$| I | $ is the length of $I$.
\end{longlist}
\end{theorem}

\subsection*{Acknowledgements} The authors would like to thank the
referee for his/her
comments which resulted in a new proof of Proposition 14.
This work was partly supported by the program Volubilis ``Action Int\'
egr\'ee
MA/06/142'' and by the Commission of the European Communities
Grant PIRSES-GA-2008-230804 within the program `Marie Curie Actions'.

\printhistory

\end{document}